\documentclass[a4paper,11pt]{article}

\usepackage{amsmath,amsfonts,amssymb,caption,graphicx}
\usepackage{ushort}
\usepackage{makeidx}
\usepackage{float}
\usepackage{accents}
\usepackage{color}
\usepackage{xcolor}
\usepackage{framed}
\usepackage{versions}
\usepackage{emptypage}
\usepackage{wrapfig}
\usepackage[T1]{fontenc}
\usepackage[utf8]{inputenc}
\usepackage{titlesec}
\usepackage{fancyhdr}
\usepackage{extramarks}
\usepackage[bookmarks]{hyperref}
\usepackage{hyperref}
\usepackage{bbm}
\usepackage{enumitem}

\usepackage{ragged2e}

\hypersetup{
    colorlinks=true,   	
    linkcolor=red,      
    citecolor = [rgb]{0 0.7 0},   	
    filecolor=magenta, 	
    urlcolor=blue
}

\newsavebox{\ideabox}
\newenvironment{numberedquote}
  {\begin{equation}
   \begin{lrbox}{\ideabox}
   \begin{minipage}{\dimexpr\columnwidth-2\leftmargini}
   \setlength{\leftmargini}{0pt}%
   \begin{quote}}
  {\end{quote}
   \end{minipage}
   \end{lrbox}\makebox[0pt]{\usebox{\ideabox}}
   \end{equation}}

\title{
\vspace*{-0.1in}
Shakespeare, Entropy and Educated Monkeys}

\author{Ioannis Kontoyiannis}

\date{\today}

\begin{document}

\maketitle

\footnotetext{Statistical Laboratory, Centre for Mathematical Sciences, 
University of Cambridge, Wilberforce Road, Cambridge CB3 0WB, 
UK. Email: \texttt{yiannis@maths.cam.ac.uk}.}

\begin{abstract}
It has often been said, correctly, that a monkey forever randomly
typing on a keyboard would eventually produce the complete
works of William Shakespeare. Almost just as often it has been pointed
out that this ``eventually'' is well beyond any conceivably
relevant time frame. We point out that an {\em educated}
monkey that still types at random but is constrained to
only write ``statistically typical'' text,
would produce any given passage in a much shorter
time. Information theory gives a very simple way to estimate
that time. For example, Shakespeare's phrase,
\texttt{Better three hours too soon than a minute too late},
from {\em The Merry Wives of Windsor},
would take the educated monkey only 73 thousand years to produce,
compared to the beyond-astronomical $2.7\times 10^{63}$ years
for the randomly typing one. Despite the obvious improvement,
it would still take the educated monkey 
an unimaginably long $10^{42,277}$
years to produce all of {\em Hamlet}.
\end{abstract}


\section{The infinite monkey theorem}

It is an elementary mathematical fact that anything which
is not impossible will certainly happen eventually:
However unlikely a potential event, if 
its probability is nonzero,
then in repeated trials it is bound to eventually occur. 
To emphasize
this point, \'{E}mile Borel wrote in 1913~\cite{borel:13}:
\begin{quote}
``Imagine that a million monkeys were trained to randomly type on 
a typewriter and 
that [...] these typing monkeys worked diligently for ten hours a day 
with a million 
typewriters [...] After a year, their volumes 
would contain exact 
copies of books of every kind and in every language held in 
the richest libraries 
in the world.''
\end{quote}

\noindent
The image of a monkey sitting at a typewriter 
quickly captured the popular imagination,
and certain adaptations of
Borel's claim
became colloquially known as the {\em Infinite monkey theorem},
with references to variations of the following statement~\cite{wiki:monkeys}
regularly appearing in popular culture up to this day:
\begin{quote}
``The infinite monkey theorem states that a monkey hitting keys 
independently and at random on a typewriter keyboard for an 
infinite amount of time will almost surely type [...]
the complete works of William Shakespeare.''
\end{quote}

\noindent
For example, 
a satirical article titled
``The Infinite-Monkey Theorem: Field Notes'' appeared
in January 2023 
in {\em The New Yorker}~\cite{perlman:23},
the BBC news item
``Monkeys Will Never Type Shakespeare, 
Study Finds'' was posted online in November 2024~\cite{BBC:24},
and in the 2025 Hollywood movie {\em Superman}, Lex Luthor 
employed hundreds of super-powered monkeys to write negative comments about 
Superman on social media.

\section{The {\em finite} monkey theorem}

Borel's fantastical thought experiment and its numerous
variants received
a lot of attention at the time, though reactions were mixed.
R.G. Collingwood in his 1938 book ``The Principles of 
Art''~\cite{collingwood:58} took issue with ``random''
works of literature. He argued that art cannot be produced by accident, 
and that a piece of intentional human writing is not the same as 
one randomly generated, even if they are identical. Interestingly, 
this point is still very much relevant today,
in relation to current debates about
the value of art created by AI agents based on neural
networks and large language models~\cite{thomson:15,nannicelli:25}.
In a different direction, the literary polymath
Jorge Luis Borges in his 1938 essay ``The Total Library''~\cite{borges:39} 
traced the history of Borel's idea from Aristotle's {\em On Generation and 
Corruption} and Cicero's {\em De Natura Deorum}, through 
the writing of Blaise Pascal and Jonathan Swift, and up to 20th 
century statements with their emblematic monkeys on typewriters.

\medskip

\noindent
{\bf Yes, but when?}
Most of the discussion following Borel's observation
and its modern incarnations centers around the
question of their practical relevance.
Kittel and Kroemer~\cite{kittel:book} estimate
the probability of a monkey randomly typing {\em Hamlet} within 
the age of the universe at about $10^{-146,316}$.
Not great. 
As the authors conclude, ``the probability of {\em Hamlet}
is therefore zero in any operational sense.''
Similarly, Styer~\cite{styer:book} 
estimates that the probability of {\em Hamlet} appearing 
in the text typed by $10^{10}$ monkeys in $10^{18}$ seconds
(about $300,000,000,000$ years),
is approximately $10^{-146,211}$. 
Again, not promising. 

There is more bad news.
The evolutionary biologist and best-selling author 
Richard Dawkins
noted~\cite{dawkins:book}
that the probability of 
an even very much more modest goal -- namely, 
the appearance in a random text of the 28-character long 
phrase \texttt{METHINKS IT IS LIKE A WEASEL}
from {\em Hamlet} -- is only 
about $10^{-40}$.
And the acclaimed
populariser of mathematics
Ian Stewart 
estimated~\cite{stewart:book} 
that, in order to have a decent chance 
to observe the 8-character long phrase
\texttt{DEAR SIR},
you need to wait until
approximately 2,821,109,907,456 characters have been
typed. 

A careful and systematic evaluation of
the time it would take, on the average,
for a given target passage to appear in 
a randomly generated text was performed
by Woodcock and Falletta in 2024~\cite{woodcock:24}. 
They seem to be the ones who coined
the term {\em Finite monkey theorem},
and they estimated that the expected number of 
keystrokes needed before we 
observe Shakespeare's corpus
is approximately ``1 followed by 7.45 million zeros.''
From this they (quite reasonably)
concluded that the Infinite monkey theorem 
belongs ``in a class of probabilistic problems or paradoxes 
[...] wherein the infinite-resource 
conclusions directly contradict those obtained when 
considering limited resources, however sizeable.''
Most relevant for our purposes is their justification
of the following intuitively obvious rule of thumb:

\smallskip

\begin{numberedquote}
\centering
On the average, 
it takes about $m^{\ell}$ keystrokes
for an arbitrary text of $\ell$ characters
from an alphabet of $m$ symbols
to be
produced by a monkey typing at random.
\label{eq:woodcock}
\end{numberedquote}

\medskip

Despite all this very clear (and very easy to obtain) theoretical evidence
pointing to the impossibility of meaningful results
appearing in any kind of a realistic time frame, some people were not 
discouraged.

\medskip

\noindent
{\bf Experiments.}
In 2003, an online experimental project was created,
{\em The Monkey Shakespeare Simulator}~\cite{barrow:05,acocella:07}.
Distributed simulations
of ``virtual monkeys'' were run by various
groups, looking for text 
that matched passages from Shakespeare.
They famously reported a 24-character hit from {\em Henry IV} after
``2,737,8502,737,850 million billion billion billion monkey-years''
in equivalent simulation time. The identified phrase
consists of the first 24 characters of the randomly 
produced text

\begin{center}
\texttt{RUMOUR. 
Open your ears; 9r"5j5\&?OWTY Z0d "B-nEoF.v}\ldots
\end{center}

Around the same time, an experiment using real monkeys 
was performed at
Plymouth University, as part of the
MediaLab Arts course; see, e.g., the
BBC report~\cite{BBC:03}. Over a period of one
month, six macaques with a keyboard 
generated five pages of what appeared nothing like
random text: The result mostly looked like

\begin{center}
\texttt{SSSSSSSSSSSSSSSSSSSSSSSSSSSSSSSSSSSSS}
\end{center}

\noindent
Eventually, the lead male began striking the keyboard with a stone, 
and other monkeys followed by urinating 
on it.  Mike Phillips, the director of the university's 
Institute of Digital Arts and Technology 
rather inevitably 
concluded that real monkeys ``are not random generators''~\cite{Wired:03}.

How about {\em educated} monkeys?

\section{Educated monkeys and their entropy}

Instead of a monkey trained to type completely randomly, 
suppose we consider an {\em educated monkey} (or, for that matter,
a {\em human})
who randomly types {\em statistically typical} text that is
grammatically and syntactically correct. 
Since this excludes all meaningless 
strings of characters, we naturally expect that the time
it takes to generate any one particular phrase 
or passage should be much shorter than with purely random typing.
So, how long would
it take before we see all of {\em Hamlet}, say?

This question is difficult to answer directly because
we do not have an accurate description, or even a good 
approximation, for the probability distribution of typical
English text. 
Luckily, information theory offers a way around this obstacle.

Suppose that, instead of its full distribution,
we know the 
{\em entropy rate} $h$~\cite{cover:book2}
of typical English
text. 
Shannon's celebrated Asymptotic Equipartition Property 
(AEP)~\cite{shannon:48,cover:book2} then tells us that
there are approximately $2^{\ell h}$ typical texts
of length $\ell$, all with essentially the same
probability of around $2^{-\ell h}$ of occurring. 
Following the same logic as in~\cite{woodcock:24},
on the average it would take
$2^{\ell h}$ keystrokes to produce a given typical text
of $\ell$ characters.

To translate this into time units,
we assume our monkey types at the same speed 
of 52 words per minute
as the average human~\cite{dhakal:18},
and we adopt
the standard convention
in text-entry research that the
average word-length (including spaces)
is five characters~\cite{mackenzie:02}.
This gives a typing speed of
$\frac{13}{3}$ characters per second
or, assuming the monkey types
continuously for 24 hours a day every
day, $1.367\times 10^8$ characters per year.

Finally, we need an estimate of the entropy
rate $h$ of English text. Estimating
$h$ has a long and rich history, 
some of which is briefly outlined
in Appendix~\ref{app}. We will use what
appears to be one of the most reliable
current estimates,
$h=0.863$ bits per character, 
from~\cite{guo:24}.

Combining all these, gives the
following rule of thumb: 

\smallskip

\begin{numberedquote}
\centering
On the average, 
it takes about
$7.3\times 10^{0.26\ell -9}$ years
for a {\em typical} English text\\
of $\ell$
characters
to be produced by an educated monkey typing at random.
\label{eq:typical}
\end{numberedquote}

\medskip

How does this compare with the randomly typing monkey's time
in~(\ref{eq:woodcock})?
In order to give them the best chance,
we assume that the randomly typing monkey uses
an alphabet of only 27 characters
that consists of all lower-case letters plus space -- ignoring
all numerals, punctuation marks, and special characters.
With $m=27$ and applying the same characters-to-years conversion
as before, the estimate in~(\ref{eq:woodcock}) becomes
$7.3\times 10^{1.43\ell -9}$ years.

Simply comparing
$10^{1.43\ell}$ 
with $10^{0.26\ell}$,
it is  obvious that the educated
monkey has an enormous advantage:
An exponential with exponent 1.43 grows extremely more
rapidly than one with exponent 0.26.
To illustrate this, we take a look at some specific
examples.

\section{The importance of education}

Before examining full-length texts like {\em Hamlet},
first we consider some shorter phrases and passages.

Perhaps the shortest poem in the English language
is Muhammad Ali's, \texttt{Me, we}. On the average,
it would take
around 4.6 seconds for the educated monkey to
produce this, compared to 
38.2 days for the randomly typing monkey.\footnote{Recall that,
for the computation
of average waiting times, all phrases and passages considered here are
transformed to plain 27-character text by removing all punctuation,
numerals and special characters, and making all letters
lower-case.}
The {\em Terminator's} iconic, 
\texttt{I'll be back}, would appear in
2.8 minutes in the educated monkey's text,
as opposed to the
40,581,179 years it would take the randomly
typing one to produce it.
And for Julius Caesar's,
\texttt{The die is cast},
we would have to wait just a bit over half an hour, 
instead of
$2.2\times 10^{13}$ years.
At this point the difference is very clear.
In fact, since the random monkey's time has exceeded
our current estimate of  
$1.4\times 10^{10}$ years for the
age of the universe~\cite{choi:20},
let us continue with just the educated
monkey for a bit.

It would take around 8.4 days for the educated monkey
to produce Yoda's,
\texttt{May the Force be with you},
but for the title of this paper we would have
to wait 183.4 years.
Franklin D. Roosevelt's famous 1933 quote,
\texttt{The only thing we have to fear is fear itself}
would appear after
3,658 years,
the Spice Girls' chorus,
\texttt{I'll tell you what I want, what I really really want}
after 
73,000 years,
and the Rolling Stones',
\texttt{I can't get no satisfaction, gonna try and I try 
and I try and I try},
after approximately
$10^9$ years.

For longer phrases,  the educated monkey's time will
also exceed the age of the universe, but let us go on
with some inspiring quotes.
Winston Churchill's,
\texttt{Success is not final, failure is not fatal. It is the courage to 
continue that counts}, would appear after
$2.7\times10^{13}$ years, and Stephen Hawking's,
\texttt{We are just an advanced breed of monkeys on a minor planet of a 
very average star. But we can understand the Universe},
after $10^{22}$ years.
For the entire poem
{\em ``Hope'' is the thing with feathers}. by Emily Dickinson,
we would need to wait
$5.5\times 10^{79}$ years,
and for the popular lullaby {\em Twinkle twinkle little star},
$8.4\times 10^{142}$ years.

By now we have exceeded the largest imaginable time scale
by many orders of magnitude:
The time it will take for even the largest supermassive black holes
to evaporate via Hawking radiation is $10^{106}$ years \cite{frautschi:82}.
So let us wrap up by getting back to our starting point. Shakespeare.
Hamlet's existential conundrum,
\texttt{To be or not to be}, would take the
educated monkey  3 hours and 4 minutes to produce,
a much shorter time than the randomly typing monkey's
$4.2\times 10^{17}$ years.
But for all of {\em Hamlet}, while the corresponding
times are vastly different -- $10^{42,277}$ 
versus $10^{232,784}$ years -- they are both
still equally entirely hopeless.

\appendix

\section{Appendix: Estimates of the entropy of English}
\label{app}

The question of estimating
the entropy rate $h$ of English text 
was already examined by 
Shannon 
in his seminal
1948 paper~\cite{shannon:48} where he introduced
information theory.
Soon after that, 
in 1951~\cite{shannon:51}, he used two methods
to estimate $h$, both of which still 
remain relevant today. The first method
was based on approximations to the frequencies
of $n$-grams of words and letters in text, and the
second used human subjects attempting 
to predict the next letter in a given sentence. 
Based on these, he reported
estimates for $h$ ranging between 
$0.6$ and $1.3$ bits per character (bpc).
Cover and King in 1978~\cite{cover:king} 
used a related but different method, where
human subjects were asked to bet 
on the next character,
to obtain estimates for $h$ between $1.29$ and
$1.9$ bpc.
And a 2019 large-scale replication of Shannon's 
guessing experiment~\cite{ren:19} produced the estimate
$h\approx 1.22$.

In a different direction, numerous statistical
estimates using ideas from data compression
have been developed, including the work reported
in~\cite{kontoyiannis:96,kasw} which found
$h\approx 0.92-2.15$ bpc based on different
texts. One of the most reliable current
statistical estimates is $h\approx 1.13$ bpc,
obtained by the Prediction by Partial Matching
(PPM) algorithm followed 
by extrapolation~\cite{takahira:16}.

Compression algorithms have also been used directly
to give upper bounds to $h$.
One of the most interesting such results is from the Hutter 
prize project~\cite{HutterPrize},
an ongoing compression contest on the 1~GB file 
\texttt{enwik9}~\cite{MahoneyTestData}.
Because this file
consists of the first $10^9$ bytes of 
a 2006 English Wikipedia XML dump in plain text 
format (including markup, digits, 
punctuation, mixed-case letters, and some UTF-8 non-ASCI
characters), the best current estimate
of 0.887~bits/byte~\cite{HutterPrizeWinners2024}
from this project
is file-specific and does not necessarily represent
typical English text.

Some of the most reliable estimates for $h$ to date 
are perhaps the ones obtained by learning-based
methods, as reported, for example, 
in~\cite{deletang:24} and~\cite{guo:24,li:25}.
In these experiments, entropy estimates are 
computed based on the compression performance 
of various different architectures and 
large language models, 
including Meta's Llama~2, 
DeepMind's Chinchilla, Mistral AI's 7B parameter model, 
OpenAI's GPT-2, and Facebook's OPT-IML. The 
$h = 0.863$ bpc estimate we used 
earlier is from~\cite{guo:24}, obtained
using the Mistral 7B model~\cite{jiang:23}.

\section*{Acknowledgements} 
I am grateful to Michalis Loulakis and Venkat Anantharam for their
useful comments on an earlier version of this document. 

{\small
\bibliographystyle{plain}

\def\cprime{$'$}

}

\end{document}